\documentstyle[a4wide,11pt]{article}

\input{amssym.def}
\input{amssym.tex}
\newcommand{\rast}{\; \vec{\ast} \; }

\newcommand{\rop}{\; \vec{\oplus} \;}
\input xypic 
\newtheorem{example}{Example}[section]
 
\newtheorem{rem}[example]{Remark} 
\newtheorem{prop}[example]{Proposition}
\newtheorem{cor}[example]{Corollary}
\newtheorem{thm}[example]{Theorem}

\newtheorem{blank}[example]{}
\newcommand{\sqdiagram}[8]{ \diagram  #1  \rto^{#2} \dto_{#4}&
#3  \dto^{#5} \\ #6    \rto_{#7}  &  #8   \enddiagram  }
\newcommand{\wt}[1]{\widetilde{#1}} 
 
\newcommand{\wh}[1]{\widehat{#1}}
\newcommand{\Ker}{\mbox{Ker}\; }
\newcommand{\sast}{_{\ast}}
 
\newcommand{\bz}{{\Bbb Z}}
\newcommand{\Cok}{\mbox{Coker} \; }
\newenvironment{proof}{\noindent {\bf Proof} }{ \hfill $\Box$

\mbox{}}
 
\begin{document}
 \title{\large\bf  On finite induced crossed modules, \\ and the
homotopy 2-type of mapping cones}      \author{Ronald Brown  and
Christopher D. Wensley\\ School of Mathematics \\ University of
Wales, Bangor \\Gwynedd LL57 1UT,  U.K.\\{\it Dedicated to the
memory of J.H.C.Whitehead}}

\maketitle

\begin{abstract} Results on the finiteness of induced crossed
modules are proved both algebraically and topologically.  Using
the Van Kampen type theorem for the fundamental crossed module,
applications are given to the 2-types of mapping cones of
classifying spaces of groups. Calculations of the cohomology
classes of some finite crossed  modules are given, using crossed
complex methods. \end{abstract}

 \section*{Introduction}

Crossed modules were introduced by J.H.C. Whitehead in
\cite{W3}. They form a part of  what can be seen as his
programme of testing the idea of extending to higher dimensions
the methods of  combinatorial group theory of the 1930's, and of
determining some of the extra structure that was necessary to
model the geometry. Other papers of Whitehead of this era show
this extension of combinatorial group theory tested in different
directions. 

In this case he was concerned with the algebraic properties
satisfied  by the boundary map $$\partial : \pi_2(X,A) \to
\pi_1(A)$$ of the second relative homotopy group, together with
the standard action on it of the fundamental group $\pi_1(A).$
This is the {\em fundamental crossed module} of the pair $(X,A).
$ In order to determine the second homotopy group of a
$CW$-complex, he formulated and proved the following theorem for
this structure: 

\noindent {\bf Theorem W}  {\em Let $X=A\cup \{e^2_{\lambda }
\}$  be obtained from the connected space $A$ by attaching
2-cells. Then the second relative homotopy group $\pi _2(X,A) $
may be described as the free crossed $\pi_1(A)$-module on the
2-cells.} 

\noindent  The proof in \cite{W1} uses transversality and knot
theory ideas from the previous papers \cite{W2,W3}. See
\cite{B0} for an exposition of this proof. Several other proofs
are available. The survey by Brown and Huebschmann \cite{BHu},
and the book edited by Hog-Angeloni, Metzler and Sieradski
\cite{HMS}, give wider applications. 

The paper of Mac Lane and Whitehead \cite{MW} uses Theorem W to
show that the 2-dimensional homotopy theory of pointed,
connected $CW$-complexes  is completely modelled by the theory
of crossed modules. This is an extra argument for regarding
crossed modules as 2-dimensional versions of groups. 

One of our aims is the explicit calculation of examples of the
crossed module \begin{equation} \pi_2(A \cup \Gamma V,A) \to
\pi_1(A) \label{cr-m-cone} \end{equation} of a mapping cone,
when $\pi_1(V), ~\pi_1(A)$ are finite. The key to this is the
2-dimensional Van Kampen Theorem (2-VKT) proved by Brown and
Higgins in \cite{BH1}. This implies a generalisation of Theorem
W, namely that the crossed module (\ref{cr-m-cone}) is  induced
from the identity crossed module $\pi_1(V) \to \pi_1(V)$ by the
morphism $\pi_1(V) \to \pi_1(A). $  

Presentations of induced crossed modules are given in
\cite{BH1}, and from these we prove a principal theorem (Theorem
\ref{finite}), that crossed modules induced from finite crossed
modules by morphisms of finite groups are finite. We also use
topological methods to show that induced crossed modules of
finite $p$-groups are also finite $p$-groups. These results give
a new range of finite crossed modules. 

Sequels  discuss crossed modules induced by a normal inclusion
\cite{BW2}, and calculations obtained using a  group theory
package \cite{BW3}.

The origin of the 2-VKT was the idea of extending to higher
dimensions the notion of the fundamental groupoid, as suggested
in  1967 in \cite{B-1}. This led to the discovery of the 
relationship of 2-dimensional groupoids to crossed modules, in
work with Spencer \cite{BS}. This relationship reinforces the
idea of  `higher dimensional group theory', and was essential
for the proof of the 2-VKT for the fundamental crossed module. 
In view of the results of Mac Lane and Whitehead,  and of
methods of classifying spaces of crossed modules by Loday
\cite{L} and Brown and Higgins \cite{BH2} (see section
\ref{topapp}),  the 2-VKT allows for the explicit computation of
some homotopy 2-types, in the form of the crossed modules which
model them. 

In some cases,  the Postnikov invariant of these 2-types can be
calculated, as the following example shows. 

\noindent {\bf Corollary \ref{post}}  Let $C_n$ denote the
cyclic group of order $n$, and let $BC_n$ denote its classifying
space.  The first Postnikov invariant of the mapping cone
$BC_{n^2}\cup \Gamma BC_n$ is a generator of a cyclic group of
order $n$, namely the cohomology group $H^3(C_n,A_n)$, where
$A_n$ is a particular cyclic $C_n$-module. 

The method used for the calculation of the cohomology class here
is also of interest, since it involves a small free crossed
resolution of the cyclic group of order $n$ in order to
construct an explicit 3-cocycle corresponding to the crossed
module.  This indicates a wider  possibility of using crossed
resolutions for explicit calculations.  It is also related to 
Whitehead's use of what he called in \cite{W1} `homotopy
systems',  and which are simply free crossed complexes.

 \section{Crossed modules and induced crossed modules}

In this section, we recall the definition of induced crossed
modules, and of results of \cite{BH1} on presentations of
induced crossed modules. We then give some basic examples of
these.  \par  Recall that a {\em crossed module} is a morphism
of groups $\mu : M \to P$ together with an action $(m,p)\mapsto
m^p$ of $P$ on $M$ satisfying the two axioms  \begin{itemize} 
\item  CM1) $\mu (m^p) = p^{-1}mp$  \item CM2) $n^{-1}mn= m
^{\mu n}$ \end{itemize}  for all $m,n \in M, p \in P.$   \par 
The category ${\cal XM}$ of crossed modules has objects all 
crossed modules with morphisms  the commutative diagrams  $$
\sqdiagram{M}{\mu}{P}{g}{f}{N}{\nu}{Q} $$   in which the
horizontal  maps are crossed modules, and the pair $g,f $
preserve the action in the sense that for all $m\in M, p\in P$
we have $g(m^p)=(gm)^{fp}.$ If $P$ is a group, then the category
${\cal XM}/P$ of crossed $P$-modules is the subcategory of
${\cal XM}$ whose objects are the crossed $P$-modules and in
which a {\em morphism } $g : M\to N$ of crossed $P$-modules is a
morphism of groups such that $g$ preserves the action
($g(m^p)=(gm)^p,$ for all $m\in M, \; p\in P$), and $\nu g= \mu.$

Standard algebraic examples of crossed modules are: \par 
\noindent (i) an inclusion of a normal subgroup, with action
given by conjugation; \par  \noindent (ii) the inner
automorphism map  $\chi : M \to \mbox{Aut}\;M,$ in which $\chi m
$ is the automorphism $n\mapsto m^{-1}nm$; \par  \noindent (iii)
the zero  map $M \to P $ where   $M$ is a $P$-module; \par
\noindent (iv) an epimorphism $M \to P $ with kernel contained
in the centre of $M$ . \par  \noindent All these yield examples
of finite crossed modules. Other finite examples may be
constructed from those above, the induced crossed modules of
this paper and its sequel \cite{BW2}, and coproducts \cite{B2}
and tensor products \cite{BL,E} of crossed $P$-modules. \par 

Further important examples of crossed modules are the free
crossed modules, referred to in the Introduction, and which are
rarely finite.  They arise algebraically in considering
identities among relations \cite{BHu}, which are non abelian
forms of syzygies.   \par   We next define {\em pullback crossed
modules}.   Let $\iota :P\to Q$ be a morphism of groups. Let
$\nu : N\to Q$ be a crossed $Q$-module. Let $\nu ' :\iota
^{\ast}N \to P$ be the pullback of $N$ by $\iota $, so that
$\iota ^{\ast}N = \{ (p,n)\in P\times N | \iota p = \nu n \}$,
and $\nu ' :(p,n) \mapsto p.$ Let $P$ act on $\iota ^{\ast}N$ by
$(p_1,n)^p = (p^{-1}p_1p, n^{\iota p}).$ The verification of the
axiom CM1) is immediate, while CM2) is proved as follows: \par
\noindent Let $(p,n),(p_1,n_1) \in \iota^{\ast}N. $ Then  \[
\begin{array}{ccl}  (p,n)^{-1}(p_1,n_1)(p,n) & =&
(p^{-1}p_1p,n^{-1}n_1n) \\ &=& (p^{-1}p_1p,n_1^{\nu n}) \\ &=&
(p^{-1}p_1p, n_1^{\iota p}) \\ &=& (p_1,n_1)^{\nu ' (p,n)}. 
\end{array} \] 

\begin{prop}  The functor $\iota ^{\ast} : {\cal XM}/Q \to {\cal
XM}/P$  has a right adjoint $\iota_{\ast}$. \end{prop} 

\begin{proof}   This follows from general considerations on Kan
extensions.   \end{proof}

\mbox{}

The universal property of induced crossed modules is the
following. Let $\mu : M \to P, \linebreak \gamma : C \to Q$ be 
crossed modules. In the diagram  $$ \diagram & M \ddlto_f
\rto^{\mu}\dto^{\bar{\iota}}& P \dto^{\iota} \\ & \iota \sast M
\rto^{\delta} \dldashed^(0.3)g|>\tip & Q \\ C \urrto_{\gamma} &&
\enddiagram $$  the pair $  \bar{\iota}, \; \iota$\ is a
morphism of crossed modules  such that for any crossed
$Q$-module  \linebreak  $\gamma : C \to Q $ and morphism of
crossed modules $f,\; \iota, $ there is a unique morphism $g :
\iota \sast M \to C$ of crossed $Q$-modules such that $g
\bar{\iota}= f.$ \par    

It is a consequence of this universal property  that if 
$M=P=F(R)$, the free group on a set $R,$ and if $w: R \to Q $ is
the restriction of $\iota $ to the set $R,$ then
$\iota_{\ast}F(R)$ is the free crossed module on $w, $ in the
sense of Whitehead \cite{W1} (see also \cite{BHu,HMS,Pr}).
Constructions of this free crossed module are given in these
papers.\par   A presentation for induced crossed modules for a
general morphism $\iota$ is given in Proposition 8 of 
\cite{BH1}. We will need two more particular results. The first
is Proposition 9 of that paper, and the second is a direct
deduction from Proposition 10. \par 

\begin{prop}   If $\iota : P\to Q$ is a surjection, and $\mu :
M\to P$ is a crossed P-module, then $\iota_{\ast}M \cong
M/[M,K],$ where $K = \mbox{{\em Ker}} \; \iota $, and $[M,K]$
denotes the subgroup of M generated by all $m^{-1}m^k$ for all
$m\in M, \; k\in K.$ \label{surj}  \end{prop} 

The following term and notation will be used frequently.  Let
$P$ be a group and let $T$ be a set. We define the {\em 
copower} $P \rast T$ to be the free product of groups $P_t, t\in
T, $ each with elements $(p,t), p \in P, $ and isomorphic to $P
$ under the map $(p,t) \mapsto p. $ If $Q$ is a group, then $P
\rast Q $ will denote the copower of $P$ with the underlying set
of the group $Q.$

\begin{prop}  If  $\iota : P\to Q $ is an injection, and $\mu :
M\to P $ is a crossed P-module, let T be a right transversal of 
$\iota P$  in Q.  Let Q act on the copower $M \rast T$ by the
rule $(m,t)^q= (m^p,u),$ where $p\in P, \; u\in T,$ and
$tq=(\iota p)u. $ Let $\delta : M \rast T  \to Q$ be defined by
$(m,t) \mapsto t^{-1}(\iota \mu m)t.$  Let S be a set of
generators of M as a group, and let $S^P= \{ x^p:x\in S,\;p\in P
\}$. Then $$\iota_{\ast}M = (M \rast T) /R$$ where R is the
normal closure in $M \rast T $  of the elements  $$\langle
(r,t),(s,u) \rangle = (r,t)^{-1}(s,u)^{-1}(r,t)(s,u)^{\delta 
(r,t)} \; \; \; \;  (r,s\in S^P, \; t,u \in T).$$ \label{inj} 
\end{prop} 

\begin{proof}  Let $N = M \rast T.$ Proposition 10 of \cite{BH1}
yields that $\iota_{\ast}M$ is the quotient of  $N$ by the
subgroup $\langle N, N \rangle$ generated by $\langle
n,n_1\rangle = n^{-1}n_1^{-1}nn_1^{\delta n}, n,n_1 \in N, $ and
which is called in \cite{BHu} the Peiffer subgroup of $N$.  Now
$N$ is generated by the set $(S^P,T) = \{ (s^p,t) : s\in S, p\in
P, t\in T\}, $ and this set is $Q$-invariant since $(s^p,t)^q =
(s^{pp'},u)$ where $u \in T, p' \in P$ satisfy $tq = (\iota
p')u.$  It follows from Proposition 3 of \cite{BHu} that
$\langle N,N \rangle$ is the normal closure of the set $ \langle
(S^P,T), (S^P,T) \rangle $ of basic Peiffer commutators.
\end{proof}  

\begin{example} The dihedral crossed module.  {\em We show how
this works out in the following case, which exhibits a number of
typical features. We let $Q$ be the dihedral group $D_n$ with
presentation $\langle x,y : x^n = y^2 =xyxy =1 \rangle$, and let
$M=P$ be the cyclic subgroup $C_2$ of order 2 generated by $y.$
Let $C_n = \{ 0,1,2,\ldots,n-1 \} $ be the cyclic group of order
$n.$ A right transversal $T$ of $C_2$ in $D_n$ is given by the
elements $x^i, \; i \in C_n .$ Hence $\iota_{\ast}C_2$  has a
presentation with generators $a_i=(y,x^i), \; i\in C_n, $  and
relations given by $a_i^2 =1, i \in C_n$, together with the
Peiffer relations. Now $\delta a_i = x^{-i}yx^i = yx^{2i}. $
Further the action is given by $(a_i)^x = a_{i+1}, \; (a_i)^y =
a_{n-i}. $ Hence $(a_i)^{\delta a_j} = a_{2j-i}$, so that the
Peiffer relations are $a_ja_ia_j= a_{2j-i}. $ It is well known
that we now have a presentation of the dihedral group $D_n$, in
which we get the standard presentation $\langle u,v : u^n= v^2 =
uvuv=1 \rangle$ by setting $u=a_0a_1, \; v= a_0, $  so that $u^i
= a_0a_i.$  Then $$\delta u = x^2, \; \delta v = y,$$  so that
$y$ acts on $\iota_{\ast}C_2$ by conjugation by $v.$ However $x$
acts by   $$u^x=u,\;  v^x= vu.$$  Note that this is consistent
with the crossed module axiom CM2) since 
$$v^{x^2}=(vu)^x=vuu=u^{-1}vu.$$  We call this crossed module
the {\em dihedral crossed module}. It follows from these
formulae that $\delta $ in the induced crossed module $\delta :
D_n \to D_n $ is an isomorphism if $n$ is odd, and has kernel
and cokernel isomorphic to $C_2$ if $n$ is even. In particular,
if $n$ is even, then by results of section \ref{topapp}, $\pi_2
(BD_n \cup \Gamma BC_2)$  can be regarded as having one
non-trivial element represented by $u^{n/2}$.  }
\label{dihedral} \end{example} 

\begin{cor} Assume $\iota : P \to Q$ is injective.   If M has a
presentation as a group with $g$ generators and $r$ relations,
the set of generators of M is P-invariant, and $n=[Q:\iota \mu
(M)],$ then $\iota_{\ast}M$ has a presentation with $gn$
generators and $rn+g^2n(n-1)$ relations.  \end{cor}

Another corollary determines induced crossed modules under some
abelian conditions. This result has useful applications.  If $M$
is an abelian group, or $P$-module, and $T$ is a set we define
the {\em copower} of $M$ with $T$, written $M \; \vec{\oplus} \;
T $, to be the sum of copies of $M$ one for each element of $T.$ 

\begin{cor} Let $\mu : M \to P $ be a crossed $P$-module and
$\iota : P \to Q $ a monomorphism of groups such that M is
abelian and $\iota  \mu (M) $ is normal in $Q$. Then
$\iota_{\ast}M$ is abelian and as a $Q$-module is just the
induced $Q$-module in the usual sense.   \label{Mabelian} 
\end{cor} 

\begin{proof} We use the result and notation of Proposition 1.3.
Note that if $u,t \in T$ and $r\in S $ then $u\delta (r,t) =
ut^{-1}(\iota  \mu r)t=(\iota  \mu m)ut^{-1}t = (\iota  \mu m)u$
for some $m \in M,$ by the normality condition. The Peiffer
commutator given in Proposition 1.3 can therefore be rewritten
as  $$(r,t)^{-1}(s,u)^{-1}(r,t)(s,u)^{\delta  (r,t)} =
(r^{-1},t) (s,u)^{-1}(r,t)(s^{m},u).  $$   Since $M$ is abelian,
$s^{m}=s.$ Thus the basic Peiffer commutators reduce to ordinary
commutators. Hence $\iota_{\ast}M $ is the copower $ M \;
\vec{\oplus}\; T,$ and this, with the given action,  is the
usual presentation of the induced $Q$-module. \end{proof}

\begin{example}  {\em  Let  $M=P=Q$ be  the infinite cyclic
group, which we write $\Bbb Z$, and let $\iota : P\to Q $ be
multiplication by 2.   Then $\iota_{\ast}M \cong {\Bbb Z} \oplus
{\Bbb Z}, $ and the action of a generator of $Q$ on
$\iota_{\ast}M $ is to switch the two copies of $\Bbb Z$. This
result could also be deduced from well known results on free
crossed modules. However, our results show that we get a similar
conclusion simply by replacing each $\Bbb Z$ in the above by for
example $C_4 $, and this fact is new.  }  \end{example}

 \section{On the finiteness of induced crossed modules} 

In this section we give an algebraic proof that a crossed module
induced from a finite crossed module by a morphism with finite
cokernel is also finite. In a later section we will prove a
slightly less general result, but  by topological methods which
will also yield results on the preservation of the Serre class
of a crossed module under the inducing process. 

\begin{thm}   Let $\mu : M \to P  $ be a crossed module and let
$\iota : P \to Q $ be a morphism of groups. Suppose that $M$  
and the index of $\iota (P) $ in $Q$ are finite. Then the
induced crossed module $\iota \sast M $ is finite.
\label{finite}  \end{thm} 

\begin{proof} Factor the morphism $\iota : P \to Q $ as $ \tau
\sigma $ where $\tau$  is injective and $\sigma $ is surjective.
Then $\iota \sast M $ is isomorphic to $\tau \sast \sigma \sast
M. $ It is immediate from Proposition \ref{surj} that if $M$ is
finite then so also is $\sigma \sast M. $ So it is enough to
assume that $\iota $ is injective, and in fact we assume it is
an inclusion. \par Let $T$ be a right transversal of $P$ in $Q.$
Then there is a function  $$(\xi , \eta) : T \times Q \to P
\times T $$ determined by the equation  $$tq = \xi (t,q) \eta
(t,q) $$  for all $t\in T, \, q \in Q. $ Let $Y  = M \rast T $
be the copower of $M$ and $T,$ and let $\delta : Y \to Q $ and
the action of $Q$ on $Y$ be as in Proposition \ref{inj}. A basic
Peiffer relation is then of the form  \begin{equation}  
(m,t)(n,u) = (n,u) (m^{\xi (t, u^{-1}(\mu m) u)},  \eta ( t,
u^{-1}\mu (m) u))   \label{bpeiff}  \end{equation}  where $m,n
\in M, \, t,u \in T. $ \par  We now assume that the finite set
$T$ has been given the total order $t_1 < t_2 < \ldots < t_l$.
An element of $Y$ may be represented as a reduced word 
\begin{equation} (m_1,u_1)(m_2,u_2) \ldots (m_e,u_e).
\label{word}  \end{equation} Such a word is said to be {\em
ordered} if $u_1 < u_2 < \ldots < u_e$ in the given order on $T.
$ \par \noindent  {\bf Claim: } {\em There is an algorithm
which, applying the Peiffer relations {\em (\ref{bpeiff})} with
$u < t $  to the reduced word {\em (\ref {word})}, yields an
ordered word. } \par  This is proved as an application of a 
purely combinatorial result in  \cite{BW1}. \par  Let $Z =
M_{t_1} \times M_{t_2} \times \ldots \times M_{t_l}$ be the
product of the sets $M_{t_i} = M\times \{ t_i \} $.  Then the
Claim implies that there is a function $\phi : Y \to Z$ such
that the quotient morphism $Y \to \iota \sast M $ factors
through $\phi. $ But $Z$ is finite. It follows that $\iota \sast
M$ is finite.   \end{proof} 

 \section{Topological applications}\label{topapp}  As explained
in  the Introduction, the fundamental crossed module functor
$\Pi_2$ assigns  a crossed module $\partial : \pi_2(X,A) \to
\pi_1(A)$ to any base pointed pair of spaces $(X,A)$. We will
use the following consequence of Theorem C of \cite{BH1}, which
is a 2-dimensional Van Kampen type theorem for this functor.
\begin{thm}  {\em (\cite{BH1}, Theorem D)}  Let $(B,V)$ be a
cofibred pair of spaces,  let $f : V \to A $ be a map, and let 
$X=A \cup_f B $. Suppose that $A,B,V$ are path-connected, and
the pair $(B,V)$ is $1$-connected. Then the pair $(X,A)$ is
$1$-connected and the diagram  $$
\sqdiagram{\pi_2(B,V)}{\delta}{\pi_1(V) } {\epsilon }{\lambda }
{\pi_2(X,A) }{\delta'}{\pi_1(A)} $$ \noindent presents
$\pi_2(X,A)$ as the crossed $\pi_1(A)$-module
$\lambda_{\ast}(\pi_2(B,V))$ induced from the crossed
$\pi_1(V)$-module $\pi_2(B,V)$ by the group morphism $\lambda :
\pi_1(V) \to \pi_1(A)$  induced by $f.$ \label{VKT}  \end{thm}
\par  

As pointed out earlier, in the case $P$ is a free group on a set
$R$,  and $\mu$ is the identity, then the induced crossed module
$\iota_{\ast}P$ is the free crossed $Q$-module on the function
$\iota |R :R \to Q $. Thus Theorem \ref{VKT} implies Whitehead's
Theorem W of the Introduction.  A considerable amount of work
has been done on this case, because of the connections with
identities among relations, and methods such as transversality
theory and ``pictures''  have proved successful (\cite{BHu,Pr}),
particularly in the homotopy theory of 2-dimensional complexes
\cite{HMS}. However, the only route so far available to the
wider geometric applications of induced crossed modules is
Theorem \ref{VKT}.  We also note that this Theorem includes the
relative Hurewicz Theorem in this dimension, on putting $A =
\Gamma V$, and $f : V \to \Gamma V$ the inclusion.   \par 

We will apply this Theorem \ref{VKT} to the {\em classifying
space of a crossed module}, as defined by Loday in \cite{L} or
Brown and Higgins in \cite{BH2}.  This classifying space is a
functor $B$ assigning to a crossed module ${\cal M}= (\mu : M
\to P)$ a pointed $CW$-space $B{\cal M}$  with the following
properties:

\noindent   \begin{blank} The homotopy groups of the classifying
space of the  crossed module $\mu : M\to P $ are given by  
$$\pi_i(B(M \to P)) \cong \left\{  \begin{array} {ll}
\mbox{Coker}\; \mu & \mbox{for $i=1$} \\      \mbox{Ker}\;   \mu
& \mbox{for $i=2$} \\             0           & \mbox{for
$i>2.$}   \end{array} \right. $$  \end{blank}

\noindent  \begin{blank}  The classifying space $B(1\to P)$ is
the usual classifying space  $BP$ of the group $P$,  and $B P$
is a subcomplex of $B(M \to P).$  Further, there is a natural
isomorphism of crossed modules $$\Pi_2(B(M\to P),BP) \cong (M
\to P).$$ \end{blank} 

\noindent  \begin{blank}  If $X$ is a reduced $CW$-complex with
$1$-skeleton $X^1,$  then there is a map  $$X \to
B(\Pi_2(X,X^1))$$ inducing an isomorphism of $\pi_1$ and
$\pi_2$.   \end{blank} 

It is in these senses that it is reasonable to say, as in the
Introduction, that crossed modules model all pointed homotopy
2-types. \par We now give two direct applications of Theorem
\ref{VKT}. 

\begin{cor}  Let  $\mu : M \to P$ be a crossed module, and let
$\iota : P\to Q$ be a morphism of groups. Let $\beta: BP \to
B(M\to P)$ be the inclusion. Consider the pushout

\begin{equation}  \diagram BP \rto^(0.33){\beta} \dto_{B \iota }
& B(M \to P)  \dto \\  BQ \rto_{\beta '} & X.  \enddiagram
\label{pushout}  \end{equation} 

Then the fundamental crossed module of the pair $(X,BQ)$ is
isomorphic  to the induced crossed module $\iota_{\ast} M \to Q
,$ and this crossed module determines the 2-type of $X.$ 
\end{cor} \begin{proof}The first statement is immediate from
Theorem  \ref{VKT}. The final statement follows from results of
\cite{BH2}, since the morphism $Q \to \pi _1(X)$ is surjective.
\end{proof}  

\begin{rem} {\em An interesting special case of the last
corollary is when $\mu : M \to P $ is an inclusion of a normal
subgroup, since then $B(M \to P)$ is of the homotopy type of
$B(P/M).$  So we have determined the 2-type of a homotopy
pushout $$\sqdiagram{BP}{Bp}{BR}{B\iota}{}{BQ}{p  '}{X}$$ in
which $p : P \to R$ is surjective. }\end{rem} 

We write $\Gamma V$ for the cone on a space $V.$

\begin{cor}  Let $\iota :P\to Q$ be a morphism of groups. Then
the fundamental crossed module   $\Pi_2(BQ \cup _{B\iota}\Gamma
BP,BQ)$  is isomorphic to the induced crossed module $\iota
\sast P \to Q.$  \end{cor} \par 

 \section{Finiteness theorems by topological methods} 

The aim of this section is to show that the property of being a
finite $p$-group is preserved by the process of induced crossed
modules. We use topological methods. \par 

An outline of the method is as follows. Suppose that $Q$ is a
finite $p$-group. To prove that $\iota \sast M$ is a finite
$p$-group, it is enough to prove that $\Ker (\iota \sast M\to
Q)$ is a finite $p$-group. But this kernel is the second
homotopy group of the space $X$ of the pushout (\ref{pushout}),
and so is also isomorphic to the second homology group of the
universal cover $\wt X $ of $X. $ In order to apply the homology
Mayer-Vietoris sequence to this universal cover, we need to show
that it may be represented as a pushout, and we need information
on the homology of the spaces determining  this pushout. So we
start with the necessary information on covering spaces. \par 

We work in the convenient category ${\cal TOP}$ of weakly
Hausdorff $k$-spaces \cite{Le}. Let $\alpha : \wt{X} \to X$ be a
map of spaces. In the examples we will use, $\alpha $ will be a
covering map. Then $\alpha $ induces a functor $$\alpha ^{\ast}:
{\cal TOP}/X \to {\cal TOP}/ \wt{X}.$$ It is known that $\alpha
$ has a right adjoint and so preserves colimits
\cite{BB1,BB2,Le}.  

For regular spaces, the pullback of a covering space in the
above category is again a covering space. These results enable
us to identify a covering space of an adjunction space as an
adjunction space obtained from the induced covering spaces. \par 

If further, $\alpha $ is a covering map, and $X$ is a
$CW$-complex, then $\wt{X}$ may be given the structure of a
$CW$-complex \cite{Ma}.   \par

We also need a special case of the basic facts on the path
components and fundamental group of induced covering maps
(\cite{B1,BHK,Ma}).   Given the following pullback   $$
\sqdiagram{\wh{A}}{}{\wt{X}}{\alpha '}{\alpha}{A}{f}{X} $$ and
points $a \in A, \; \tilde{x} \in \wt{X}$ such that $fa=\alpha
\tilde{x},$ in which $\alpha $ is a universal covering map and
$X,A,\wt{X}$ are path connected, then there is a sequence  

\begin{equation}  1\to \pi_1(\wh{A},(a,\tilde{x})) \to  \pi_1
(A,a) \stackrel{f_{\ast}}{\to} \pi_1 (X,fa) \to \pi_0 (\wh{A})
\to 1. \label{covexact}  \end{equation} 

This sequence is exact in the sense of sequences arising from
fibrations of groupoids  \cite{B1}, which involves an operation
of the fundamental group $\pi_1(X,fa)$ on the set $\pi_0
(\wh{A})$ of path components of $\wh{A}$. It follows that  the
fundamental group of $\wh{A}$ is isomorphic to Ker  $f_{\ast}$,
and that $\pi_0(\wh{A})$ is bijective with the set of cosets
$(\pi_1(X,fa))/(f_{\ast}\pi_1(A,a)). $ It is also clear that the
covering $\wh{A}\to A$ is regular and that all the components of
$\wh{A}$ are homeomorphic. \par  

Let $\cal M$ be the crossed module $\mu : M\to P$ and let
$\iota: P\to Q$ be a morphism of groups.  Let
$$X=BQ\cup_{B\iota}B{\cal M} $$ as in diagram (\ref{pushout}). 
Let $\alpha  : \wt{X} \to X $ be the universal covering map, and
let $\wh{BQ}, \; \wh{B{\cal M}}, \; \wh{BP}$ be the pullbacks of
$\wt{X}$ under the maps $BQ\to X,\; B{\cal M} \to X, \; BP \to
X.$ Then we may write 

\begin{equation} \wt{X} \cong \wh{BQ} \cup_{\wh{B\iota }}
\wh{B{\cal M}}, \label{po-xtilde} \end{equation} 

\noindent by the results of section \ref{topapp}. \par  

>From the exact sequence (\ref{covexact}) we obtain the following
exact sequences, in which  $\pi_1 X \cong Q \ast_P (P/\mu M)$:  
$$ \begin{array}{ccccccccc} 1\to & \pi_1(\wh{BQ}) & \to & Q &
\to & Q \ast_P (P/\mu M) &\to & 1, & \\ 1\to & \pi_1(\wh{B{\cal
M}}) & \to & P/ \mu M & \to & Q \ast_P (P/\mu M) & \to &
\pi_0\wh{B{\cal M}}& \to 1, \\ 1\to & \pi_1(\wh{BP}) & \to & P &
\to & Q \ast_P (P/\mu M) & \to & \pi_0\wh{BP} & \to 1.
\end{array} $$ 

\begin{prop}  Under the above situation, let the groups
$\pi_1(\wh{BP}) ,   \pi_1(\wh{B{\cal M}})$, $ \pi_1(\wh{BQ})$ be
denoted by  $P', M', Q'$ respectively, and let $B{\cal M}'$
denote a component of $\wh{B{\cal M}}$. Then there is an exact
sequence  $$ H_2(P') \rop \pi_0\wh{BP} \to (H_2(B{\cal M}')\rop
\pi_0 \wh{B{\cal M}}) \oplus H_2(Q')  \to \pi_2(X) \to $$  $$
\to  H_1(P') \rop \pi_0\wh{BP} \to ( H_1(M') \rop \pi_0
\wh{B{\cal M}}) \oplus H_1(Q')  \to 0. $$ \label{longseq}
\end{prop} 

\begin{proof}  This is immediate from the Mayer-Vietoris
sequence for the pushout  (\ref{po-xtilde}) and the fact that
$H_2(\wt{X}) \cong \pi_2(X). $   \end{proof} 

\begin{cor}  If $\iota : P\to Q$  is the inclusion of a normal
subgroup, and $X= BQ \cup_{B\iota}\Gamma BP, $ then $\pi_2(X)$
is isomorphic to $H_1(P) \otimes I(Q/ P) , $ where $I(G)$
denotes the augmentation ideal of a group $G.$   \end{cor}  

Note that this agrees with the result of Corollary 1.8 of
\cite{BW2}, in which the induced crossed module itself is
computed {\em via} the use of coproducts of crossed $P$-modules. 

\begin{cor} Let $\mu : M \to P $ be a crossed module  and let
$\iota : P \to Q $ be a morphism of groups. If M, P and Q are 
finite {\em p}-groups, then so also is the induced crossed
module $\iota \sast M $.   \end{cor}

\begin{proof}  It is standard that the (reduced) homology groups
of a finite  $p$-group are finite $p$-groups.  The same applies
to the  reduced homology of the classifying space of a crossed
module  of finite $p$-groups.  The latter may be proved using
the  spectral sequence of a covering, and Serre ${\cal C}$
theory, as in Chapters IX and X of \cite{H}. In the present
case, we need information only on $H_2(B(M\to P)),$ and some of
its connected covering spaces, and this may be deduced from the
exact sequence due to Hopf   $$ H_3K \to H_3G \to
(\pi_2K)\otimes _{\bz G} \bz \to H_2K \to H_2G\to 0 $$  for any
connected space $K$ with fundamental group $G$ (see for example
Exercise 6 on p.175 of \cite{Br}). Proposition \ref{longseq}
shows that   $\mbox{Ker} \; (\iota_{\ast}(M) \to Q)\cong
\pi_2(X)$  is a finite $p$-group. Since $Q$ is a finite
$p$-group,  it follows that $\iota \sast M$ is a finite
$p$-group.   \end{proof} 

\mbox{}

Note that these methods extend also to results on the Serre
class of an induced crossed module, which we leave the reader to
formulate.

 \section{Cohomology classes}

Recall \cite{Hu,Br} that if $G$ is a group and $A$ is a
$G$-module, then  elements of $H^3(G,A)$ may be represented by
equivalence classes of {\em  crossed sequences} 

\begin{equation} 0 \to A \to M \stackrel{\mu}{\longrightarrow }
P \to G \to 1,  \label{cross-seq}  \end{equation} 

\noindent  namely exact sequences as above such that $\mu :M \to
P $ is a crossed module. The equivalence relation between such
crossed sequences is generated by the {\em basic equivalences},
namely the existence of a commutative diagram of morphisms of
groups as follows  $$ \diagram  0 \rto & A \dto^1 \rto & M
\rto^{\mu} \dto^f & P \dto^g \rto & G \dto^1  \rto &1 \\  0 \rto
& A  \rto & M' \rto_{\mu '} & P' \rto & G  \rto & 1\enddiagram 
$$ such that $f,g$ form a morphism of crossed modules. Such a
diagram is called a {\em  morphism}  of crossed sequences. \par
The zero cohomology class is represented by the crossed sequence
$$ 0\to A \stackrel{1}{\longrightarrow}     A
\stackrel{0}{\longrightarrow} G   \stackrel{1}{\longrightarrow}
G \to 1, $$ which we sometimes abbreviate to  $$ A
\stackrel{0}{\longrightarrow} G. $$ In a similar spirit, we say
that a crossed module  $\mu :M \to P$ represents a cohomology
class,  namely an element of $H^3( \Cok \mu, \Ker \mu ). $   

\begin{example}  {\em Let $C_{n^2}$ denote the cyclic group of
order $n^2,$ written multiplicatively, with generator $u.$ Let 
$\gamma_n: C_{n^2} \to C_{n^2}$ be given by $u\mapsto u^n$. This
defines a crossed module, with trivial operations. This crossed
module represents the trivial cohomology class in
$H^3(C_n,C_n)$, in view of the morphism of crossed sequences  $$
\diagram 0 \rto & C_n \rto^1 \dto^1 & C_n \rto^0 \dto^{\lambda }
& C_n \rto^1 \dto^{\lambda} & C_n  \dto^1 \rto & 0  \\  0 \rto &
C_n \rto & C_{n^2} \rto_{\gamma _n} & C_{n^2} \rto & C_n \rto &
0   \enddiagram $$ where if  $t$ is the generator of the top
$C_n, $  then $\lambda (t) = u^n. $ } \end{example} 

\begin{example}  {\em We show that  the dihedral crossed module
$\partial : D_n\to D_n$  represents the   trivial cohomology
class. This is clear for $n$ odd, since then $\partial$ is an
isomorphism.  For $n$ even, we simply construct  a morphism of
crossed sequences as in the following diagram $$ \diagram 0
\rto& C_2 \rto^{1}\dto^{\cong} & C_2 \dto^{f_2} \rto^0 & C_2
\dto^{f_1} \rto^{1} & C_2 \dto^{\cong} \rto & 0 \\ 0 \rto & C_2
\rto & D_n \rto_{\partial} & D_n \rto & C_2 \rto &0 \enddiagram
$$ where if $t$ denotes the  non trivial element of $C_2$ then  
$f_1(t)=x, f_2(t)=u^{n/2}.$  Just for interest, we leave it to
the reader to prove that there is no morphism in the other
direction between these crossed sequences. }  \end{example} \par 

A crossed module $\mu : M \to P$ determines a cohomology class 
$$ k_{M \to P}\in H^3(\mbox{Coker}\; \mu, \Ker \mu ). $$  If $X$
is  a connected, pointed $CW$-complex with 1-skeleton $X^1,$
then the class $$k_X^3 \in H^3(\pi_1X, \pi_2X)$$  of the crossed
module $\Pi_2(X,X^1)$   is called the {\em first Postnikov
invariant} of $X.$ This class is also represented by
$\Pi_2(X,A)$ for any connected subcomplex $A$ of $X$ such that
$(X,A) $ is 1-connected and $\pi_2(A) = 0.$ It may be quite
difficult to determine this Postnikov invariant from a
presentation of this last crossed module, and even the meaning
of the word ``determine'' in this case is not so clear. There
are practical advantages in working directly with the crossed
module, since it is an algebraic object, and so it, or families
of such objects, may be manipulated in many convenient and
useful ways. Thus the advantages of crossed modules over the
corresponding 3-cocycles are analogous to some of the advantages
of homology groups over Betti numbers and torsion coefficients.

However, in work with crossed modules,  and in applications to
homotopy theory, information on the corresponding cohomology
classes, such as their non triviality, or their order, is also
of interest. The aim of this section is to give background to
such a determination, and to give two example of finite crossed
modules representing non trivial elements of the corresponding
cohomology groups.

The following general problem remains. If $G,A$ are finite,
where $A$ is a $G$-module, how can one characterise the subset
of $H^3(G,A)$ of elements represented by finite crossed modules?
This subset is a subgroup, since the addition may be defined by
a sum of crossed sequences, of the Baer type. (An exposition of
this is given by Danas in \cite{D}.)  It might always be the
whole group.

The natural context in which to show how a crossed sequence
gives rise to a 3-cocycle is not the traditional chain complexes
with operators but that of crossed complexes \cite{Hu}. We
explain how this works here.   For more information on the
relations between crossed complexes and the traditional chain
complexes with operators, see \cite{BH3}. 

Recall that a free crossed resolution of the group $G $ is a
free aspherical crossed complex $F \sast $ together with an
epimorphism $\phi : F_1 \to G $ with kernel $\delta_2 (F_2).$ 

\begin{example}  {\em  The cyclic group $C_n$ of order $n$ is
written  multiplicatively, with generator $t.$  We give for it 
a free crossed resolution $F \sast$ as follows. Set $F_1 =
C_{\infty}$, with generator written $w,$ and for $r \ge 2$ set $
F_r = (C_{\infty})^n. $ Here for $r \ge 2, \; F_r$ is regarded
as the free $C_n$-module on one generator $w_0, $ and we set
$w_i=(w_0)^{t^i}$. The morphism $\phi : C_{\infty}\to C_n$ sends
$w$ to $t,$ and the operation of $F_1$ on $F_r$ for $r \ge 2$ is
via $\phi. $   The boundaries are given by 

\begin{enumerate}  \item $\delta_2(w_i) = w ^n,$  \item for $r$
odd, $\delta_r(w_i)=w_iw_{i+1}^{-1}$,  \item for $r$ even and
greater than 2,  $\delta _r(w_i)=w_0w_1\dots w_{n-1}.$   
\end{enumerate}

Previous calculations show that $\delta_2$ is the free crossed
$C_{\infty}$-module on the element  $w^n \in C_{\infty}$. Thus 
$F\sast $ is a free crossed complex. It is easily checked to be
aspherical, and so is, with $\phi$,  a crossed resolution of
$C_n. $  }  \end{example} 

Let $A$ be a $G$-module.   Let $C(G,A,3) $ denote the crossed
complex $C$ which is $G$ in dimension 1, $A$ in dimension 3,
with the given action of $G$ on $A$, and which is 0 elsewhere, 
as in the following diagram  $$ \diagram \cdots\rto & 0 \rto & A
\rto & 0 \rto & G. \enddiagram  $$ Let $(F \sast , \phi ) $ be a
free crossed resolution of $G.$ It follows from the discussions
in \cite{BH3,BH2} that a 3-cocycle of $G$ with coefficients in
$A$ can be represented as a morphism of crossed complexes $f : F
\sast \to  C(G,A,3)$ over $\phi. $ This cocycle is a coboundary
if there is an operator morphism $l : F_2 \to A $ over $\phi :
F_1 \to G $ such that $l \delta _3 = f_3. $   $$ \diagram F_4
\dto  \rto^{\delta_4} & F_3 \dto_{f_3} \rto^{\delta_3} & F_2 
\dlto_l \dto \rto^{\delta_2} & F_1  \dto^{\phi}  \\ 0 \rto & A
\rto & 0 \rto  &  G   \enddiagram $$

To construct a 3-cocycle on $F \sast $ from the crossed sequence
(\ref{cross-seq}), first construct a morphism of crossed
complexes as in the diagram 

\begin{equation} \diagram  F_4 \rto \dto  &  F_3 \rto \dto^{f_3}
&  F_2 \rto \dto^{f_2}  &    F_1 \rto^{\phi} \dto^{f_1}  &  G
\dto^1  \\ 0 \rto  &  A \rto  &  M \rto_{\mu}  &  P \rto_{\psi} 
&  G   \enddiagram \label{cocyc} \end{equation}

\noindent   using the freeness of $F \sast $ and the exactness
of the bottom row.   Then compose this with the morphism of
crossed sequences $$ \diagram  0 \rto \dto  &  A \dto^1 \rto   
&     M \rto^{\mu} \dto  &  P \rto^{\psi} \dto^{\psi}  &  G
\dto^1 \\  0 \rto  &  A \rto  &  0 \rto  &  G \rto_1  &  G
\enddiagram  $$ Hence it is reasonable to say that the morphism
$f_3$ of diagram (\ref{cocyc}) is a 3-cocycle corresponding  to
the crossed sequence. \par 

We now use these methods in an  example. 

\begin{thm}  Let $n \ge 2,$ and let $\iota : C_n \to C_{n^2}$
denote the injection sending a generator $t$ of $C_n$ to $u^n,$
where $u$ denotes a generator of $C_{n^2.}$ Let $A_n$ denote the
$C_n$-module which is the  kernel of the induced crossed module
$\partial : \iota \sast  C_n \to C_{n^2} $. Then $H^3(C_n,A_n)$ 
is cyclic of order $n$ and has as generator the class of this
induced crossed module.  \label{class}  \end{thm}  

\begin{proof}  Write $\cal N$ for  the induced crossed module of
the theorem.  By Corollary \ref{Mabelian} the abelian group 
$\iota \sast C_n $ is the product $V=(C_n)^n$.  As a
$C_n$-module it is cyclic, with generator $v,$ say.  Write
$v_i=v^{t^i}, \, i=0,1,\ldots,n-1$.  Then each $v_i$ is a
generator of a $C_n$ factor of $V$.  The kernel $A_n$ of $\cal
N$ is a cyclic $C_n$-module on the generator $a=v_0v_1^{-1}.$ 
Write $a_i= a^{t^i}=v_iv_{i+1}^{-1}.$   As an abelian group,
$A_n$ has generators $a_0,a_1, \ldots, a_{n-1}$ with relations
$a_i^n=1, \; a_0a_1\ldots a_{n-1}=1. $          \par 

We define a morphism $f\sast : F \sast \to \cal N$ as follows.

\begin{enumerate}  \item $f_1$ maps $w$ to $u, $  \item $f_2 $
maps the module  generator $w_0$ of $F_2$ to $v=v_0.$  \item
$f_3$ maps the module generator $w_0$ of $F_3$ to $a_0.$  
\end{enumerate}   \begin{equation} \diagram (C_{\infty})^n
\dto_0 \rto^{\delta_4} & (C_{\infty})^n \dto_{f_3}
\rto^{\delta_3} & (C_{\infty})^n \dlto_l \dto^{f_2}
\rto^{\delta_2} & C_{\infty} \dto^{f_1} \rto & C_n  \dto^1 \\ 0
\rto & A_n \rto & (C_n)^n  \rto_(0.6){\nu_n } & C_{n^2} \rto &
C_n   \enddiagram \label{A-coc}  \end{equation} 

\mbox{}

The operator morphisms $f_r$ over $f_1$ are  defined completely 
by these conditions. \par 

The group of operator morphisms $g : (C_{\infty})^n \to A_n$
over $f_1$ may be identified with $A_n$ under $g\mapsto g(w_0).$
Under this identification, the boundaries $\delta _4 , \delta_3$
are transformed respectively to 0 and to $a_i \mapsto
a_i(a_i^t)^{-1}.$ So the 3-dimensional cohomology group is the
group $A_n $ with $a_i$ identified with $a_{i+1},
i=0,\ldots,n-1. $ This cohomology group is therefore isomorphic
to $C_n, $ and a generator is the class of the above cocycle
$f_3. $   \end{proof} 

\begin{cor}  The mapping cone $X = BC_{n^2} \cup_{B\iota} \Gamma
BC_n $ satisfies  $\pi _1X = C_n,$ and $\pi_2 X $ is the
$C_n$-module   $A_n$ of Theorem \ref{class}. The first Postnikov
invariant of $X$  is a generator of the cohomology group
$H^3(\pi_1 X , \pi_2X)$, which is a  cyclic group of order $n$. 
\label{post} \end{cor} 

The following is another example of a determination of a non
trivial cohomology class by a crossed module.  The method of
proof is similar to that of Theorem \ref{class}, and is left to
the reader. 

\begin{example}  {\em  Let $n$ be even. Let $C_n'$ denote the
$C_n$-module  which is $C_n$ as an abelian group but in which
the generator $t$ of the group  $C_n$ acts on the generator $t'$
of $C_n ' $ by sending it to its inverse. For $n=2,$ this gives
the trivial module. Then  $H^3(C_n,C_n') \cong C_2 $ and a
generator of this group is represented by the crossed module 
$\nu_n  : C_n \times C_n \to C_{n^2}, $ with generators
$t_0,t_1, u$ say,  and where $\nu_n  t_0 = \nu_n  t_1 =  u^n. $
Here $u \in C_{n^2}$ operates by switching $t_0, t_1. $ However,
it is not clear if this crossed module can be an induced crossed
module for $n>2$.}  \end{example}

\end{document}